\documentclass{article}[11pt]

\usepackage{amssymb,amsmath,enumerate}
\usepackage{a4,latexsym,parskip}
\usepackage{hyperref} 
\hypersetup{pdfauthor=Matthias Sch"utt}
\hypersetup{pdftitle=On the modularity of three Calabi-Yau threefolds with bad reduction at 11}

\newcommand{\cohom}[4][]{\mathrm{H}_{#1}^{#2}(#3,#4)}
\newcommand{\homology}[3]{\mathrm{H}_{#1}(#2,#3)}

\newcommand{\N} [1][] {\mathbb{N}_{#1}}

\newcommand{\Z}[1][]{\mathbb{Z}^{#1}}
\begin{document}
\setlength{\unitlength}{1cm}

\begin{center}
\LARGE{\textbf{On the modularity of three Calabi-Yau threefolds with bad reduction at 11}}
\end{center}

\vspace{0.7cm}

\begin{center}
\large{Matthias Sch\"utt}
\end{center}

\vspace{0.7cm}
\begin{small}
Abstract. This paper investigates the modularity of three non-rigid Calabi-Yau threefolds with bad reduction at 11. They are constructed as fibre products of rational elliptic surfaces, involving the modular elliptic surface of level 5. Their middle $\ell$-adic cohomology groups are shown to split into two-dimensional pieces, all but one of which can be interpreted in terms of elliptic curves. The remaining pieces are associated to newforms of weight 4 and level 22 or 55, respectively. For this purpose, we develop a method by Serre to compare the corresponding two-dimensional 2-adic Galois representations with uneven trace. Eventually this method is also applied to a self fibre product of the Hesse-pencil, relating it to a newform of weight 4 and level 27.

\end{small}

\section{Introduction}

The modularity of rigid Calabi-Yau threefolds over $\mathbb{Q}$ has recently been established for a huge %%@
class of manifolds (cf. \cite{DM}). However, the number of explicit examples of modular (rigid and %%@
non-rigid) Calabi-Yau threefolds is still quite small (cf. \cite{Y}). As a consequence, only a few %%@
primes of bad reduction, combining to the level of the associated newform, have appeared in those %%@
examples which cannot be derived from others by a change of coordinates. To my knowledge these %%@
primes are 2,3,5,7,17,73. Hence I found it worth investigating the modularity of three nodal %%@
Calabi-Yau threefolds which turn out to have bad reduction at 11. 

Furthermore these varieties, which arise as fibre products of two elliptic surfaces, are not rigid %%@
(i.e. $\mathrm{H}^{1,2}\neq 0$). Nevertheless we can apply the ideas of \cite{HV} to understand %%@
that all but a two-dimensional piece $U$ of $\mathrm{H}^3$ comes from some fibres which are %%@
interpreted in terms of elliptic curves over $\mathbb{Q}$. As these have been completely classified by %%@
their conductor in view of modularity, we are able to compute the action of Frobenius on $U$ %%@
essentially by counting points via the Lefschetz fixed point formula. Then we will identify its %%@
$L$-series (up to finitely many Euler factors) with the Mellin transform of a newform of weight 4 %%@
and level 22 and 55, respectively. 

Another new aspect in this paper is that we cannot use the well-known method of Livn\'e \cite[Thm. %%@
4.3]{L} to compare the corresponding 2-adic Galois representations as their traces are not even. %%@
Therefore we will work out an approach of Serre to this situation, as sketched in \cite{Se1}, which %%@
-as far as I know- has been applied only once before in this context by C. Schoen \cite{Sc1}. We %%@
will also show how this method together with Livn\'e's theorem eventually gives rise to an almost %%@
general criterion for comparing two-dimensional 2-adic Galois representations and use this in the %%@
very last section to prove the modularity of another non-rigid Calabi-Yau threefold, this time %%@
associated to a newform of level 27.

\textbf{Acknowledgement:} I would like to thank K. Hulek and C. Schoen for their great %%@
encouragement and J. Jones for valuable computations. 

This work was supported by the DFG-Schwerpunkt "Globale Methoden in der komplexen Geometrie".

\section{The construction}

We start by recalling the result of C. Schoen \cite{Sc2} that the fibre product $W$ of two %%@
relatively minimal, regular, rational elliptic surfaces $Y, Y'$ with sections over $\mathbb{P}^{1}$ has %%@
trivial canonical sheaf and finite fundamental group. Assuming all fibres over the set of common %%@
cusp $S''=S\cap S'$ (where $S,S'$ denote the cusps of $Y,Y'$, respectively) to be semi-stable, %%@
allows one to build a small resolution $\hat{W}$ of the singularities (ordinary double points). %%@
This construction implies $\hat{W}$ to have trivial canonical bundle and vanishing cohomology %%@
groups $\cohom{1}{\hat{W}}{\mathcal{O}_{\hat{W}}}$ and %%@
$\cohom{2}{\hat{W}}{\mathcal{O}_{\hat{W}}}$. Subject to some restrictions on the singular fibres %%@
which guarantee projectivity, $\hat{W}$ is thereby seen to be a Calabi-Yau threefold (cf. %%@
\cite[Appendix]{Y} and \cite{S} for modular examples of this kind). 

In this paper we will deal with the same construction as above with the exception of dropping this %%@
last restriction; this means, we allow the factors of the fibres containing the singularities of %%@
$W$ to be irreducible (type $I_1$ in Kodaira's notation). A consequence is that there might not %%@
exist a projective small resolution. Hence, we will consider a mixed resolution $\tilde{W}$ which %%@
blows up all nodes of fibres involving an irreducible factor in the usual way ("big") while %%@
resolving all others "small" by blowing up components of the fibres. The canonical divisor of the %%@
mixed resolution $\tilde{W}$ therefore consists of the exceptional quadrics of the big %%@
resolutions, i.e. the nodal variety, where only the (projective) small resolutions are taken into account, is Calabi-Yau. 

In order to understand the third cohomology group of $\tilde{W}$ recall the isomorphism %%@
$\mathrm{H}^{1,2}(\tilde{W})\cong\mathrm{H}^{1,2}(\hat{W})$ from Schoen's article \cite[p. %%@
189]{Sc2}. Working with an analytical small resolution $\hat{W}$, we follow Schoen's arguments to %%@
deduce

\begin{eqnarray}
h^{1,2}(\hat{W})=1+\text{rk Pic}(W_{\eta})-\# S''+\sum_{s\in S\cup S'\setminus S''}(b(s)b'(s)-1)
\end{eqnarray}

where $b(s)$ (or $b'(s)$) denotes the number of components of the fibre of $Y$ (or $Y'$) above %%@
$s.$ This equation essentially comes from comparing two expressions for the Euler number of %%@
$\hat{W}$: On the one hand, the fibration implies that it simply equals twice the number of nodes; %%@
on the other hand, we can speak about a Hodge decomposition of $\mathrm{H}^{\ast}(\hat{W})$, since %%@
$\hat{W}$ is a Moishezon manifold. Hence, its Euler number is given by twice the difference %%@
$h^{1,1}(\hat{W})-h^{1,2}(\hat{W})$. Since $h^{1,1}$ may be identified with the rank of the Picard %%@
group by virtue of the exponential sequence, a three-dimensional analogue of Zariski's Lemma (cf. %%@
\cite[p. 190]{Sc2}) immediately leads to equation (1), saying that Pic$(\hat{W})$ is generated by %%@
the Picard group of the generic fibre, the 0-section and the components of the fibres which do not %%@
meet the 0-section. 

The three examples to be considered in the next section are united by the fact that the first %%@
three summands on the right hand-side of equation (1) add up to 0. Then all of %%@
$h^{1,2}(\tilde{W})$ comes from fibres of the shape $E \times I_n$ with $E$ an elliptic curve and %%@
$n\geq 2$. This enables us to apply an idea of K. Hulek and H. Verrill, developed in \cite{HV}, to %%@
compute the action of Frobenius on these parts of $\mathrm{H}^3(\tilde{W})$, leaving only the %%@
two-dimensional orthogonal complement open.

We shall now outline the structure of this method. For this purpose let $E\times \mathbb{P}^{1}$ be one %%@
component of such a fibre of type $E\times \mathbb{P}^{1}$. Its natural embedding $\iota$ in %%@
$\tilde{W}$ induces a homomorphism

\[
\iota_{\ast}: \homology{3}{E\times\mathbb{P}^{1}}{\Z}\rightarrow\homology{3}{\tilde{W}}{\Z}.
\]

We are able to apply the arguments of \cite[p.24~ff.]{HV} to this situation and realize that the %%@
images of all such components not meeting the 0-section are linearly independent in %%@
$\homology{3}{\tilde{W}}{\Z}$ by intersection theory. Hence they span a subspace %%@
$V\subset\homology{3}{\tilde{W}}{\Z}$ of dimension $2~h^{1,2}$. We turn now to $\ell$-adic %%@
cohomology and start with the natural Galois representation $\varrho$ of $\tilde{W}$
\[
\varrho: \text{Gal}(\bar{\mathbb{Q}}/\mathbb{Q}) \rightarrow %%@
\text{Aut}(\cohom[\text{\'et}]{3}{\tilde{W}}{\mathbb{Q}_{\ell}})
\cong \text{GL}_{h^3}(\mathbb{Q}_{\ell}).
\]

By the above considerations, the induced homomorphisms $\iota^{\ast}$ for those components which %%@
do not meet the 0-section provide a short exact seqence of $\ell$-adic Galois representations
\begin{eqnarray}
0\rightarrow U \rightarrow \cohom[\text{\'et}]{3}{\tilde{W}}{\mathbb{Q}_{\ell}} \rightarrow V \rightarrow 0
\end{eqnarray}

with $U$ two-dimensional. We are going to prove the modularity of $U$, but before getting into %%@
detail we shall introduce the elliptic surfaces to be considered.

\section{The fibre products}

In this section we are going to consider three fibre products of rational elliptic surfaces with %%@
sections. These fibre products will have one factor in common which is the modular elliptic %%@
surface of level 5, $S_1(5)$. By Beauville's classification \cite{B}, a birational model is given %%@
by the hypersurface

\[
x(x-z)(y-z)+tyz(x-y)=0
\]

in $\mathbb{P}^{1}\times\mathbb{P}^{2}.$ Resolving the singularities at $(0,(0:1:0)), (0,(1:1:1)), %%@
(\infty,(1:0:0)),$ and $(\infty,(0:0:1))$, gives $S_1(5)$ the structure of an elliptic surface %%@
with the following singular fibres

\vspace{0.3cm}
\begin{center}
\begin{tabular}{c|c|c}
cusp & type of fibre & coordinates of nodes  \cr
\hline
$\infty$ & $I_5$ & rational\cr
$0$ & $I_5$ & rational\cr
$\frac{11+5\sqrt{5}}{2}$ & $I_1$ & irrational; in $\mathbb{Q}(\sqrt{5})$\cr
$\frac{11-5\sqrt{5}}{2}$ & $I_1$ & irrational; in $\mathbb{Q}(\sqrt{5})$\cr
\end{tabular}
\end{center}
\vspace{0.3cm}

Note that the six surfaces classified by Beauville have the property of being \emph{extremal} in %%@
common, i.e. the rank of the Picard group of the generic fibre is one. In view of equation (1) it %%@
is worth noticing that the Picard group of the generic fibre $\tilde{W}_{\eta}$ of a fibre product %%@
is the direct sum of the Picard groups of the factors, Pic $Y_{\eta}$ and\linebreak
Pic $Y'_{\eta}$, if $Y$ and $Y'$ are not isogenous (which is in particular the case if $S\neq %%@
S'$). As we want to look at examples where the right hand-side of equation (1) reduces to the last %%@
summand, our first aim is to build a fibre product of $S_1(5)$ with itself with three common %%@
cusps.

\subsection{A self fibre product of $S_1(5)$: $\tilde{W}_1$}

Our first variety $\tilde{W}_1$ may be realized by modifying the natural projection of $S_1(5)$ %%@
onto $\mathbb{P}^{1}$ by an automorphism $\pi$ which permutes exactly three of the four cusps. However, %%@
we need that $\pi$ is defined over $\mathbb{Q}$. This leaves -up to an isomorphism which correspond to the %%@
twist by $t \mapsto -\frac{1}{t}$- only one possibility, namely

\[
\pi:  t \mapsto 11-t.
\]

Then we consider the twisted fibre product $(S_1(5), pr)\times_{\mathbb{P}^{1}}(S_1(5),\pi\circ pr).$ %%@
The 25 ordinary double points of the fibre above $\infty$ allow to produce projective small %%@
resolutions by successive blow-ups of components of the fibre. Meanwhile the two remaining nodes %%@
in the $I_1 \times I_1$-fibres are blown up in the usual (big) way. Let $\tilde{W}_1$ denote the %%@
resulting mixed resolution. In order to determine the bad primes of $\tilde{W}_1$ we realize that %%@
there are two ways for $\tilde{W}_1$ to obtain bad reduction, since the reduction of $S_1(5)$ is %%@
everywhere good. On the one hand degenerations of fibres can produce worse singularities such that %%@
our mixed resolution is not sufficient. In this example this happens only in characteristic 5 %%@
where the two $I_1$ fibres degenerate to one of type $II$ in $\mathbb{F}_{5}$. Indeed, the resulting %%@
singularity of type $II \times II$ requires four usual blow-ups instead of the two used. On the %%@
other hand, congruences of cusps in $\mathbb{F}_{p}$ will in general cause additional double points. Here, %%@
the two non-common cusps 0 and 11 come together modulo 11, leading to 25 additional nodes. Hence %%@
the bad primes of $\tilde{W}_1$ are 5 and 11.

We will continue the investigation of the variety $\tilde{W}_1$ in the sixth section. Meanwhile we will now %%@
introduce two other fibre products which have $S_1(5)$ as one factor and an elliptic surface $Y'$ %%@
with five semi-stable singular fibres as second factor. Then Zariski's Lemma implies rk Pic %%@
$Y'_{\eta}=2$. Thus, if four of the cusps of $Y'$ coincide with those of $S_1(5)$, their fibre %%@
product will have the right hand-side of equation (1) again reduce to the last summand, as %%@
required.

\subsection{The fibre product $\tilde{W}_2$}

We consider the elliptic surface $Y'$ arising from the following pencil of cubics:

\[
(x+y+z)(\frac{11}{8}xy+\frac{11}{8}yz+\frac{125}{88}zx)-(t+\frac{125}{88})xyz.
\]

After resolving the three singularities in the fibre above $\infty$, we find that the singular %%@
fibres have the following structure 

\vspace{0.3cm}
\begin{center}
\begin{tabular}{c|c|c}
cusp & type of fibre & coordinates of nodes  \cr
\hline
$\infty$ & $I_6$ & rational\cr
$-\frac{125}{88}$ & $I_2$ & irrational; in $\mathbb{Q}(\sqrt{-5\cdot 359})$\cr
$0$ & $I_2$ & irrational; in $\mathbb{Q}(\sqrt{5})$\cr
$\frac{11+5\sqrt{5}}{2}$ & $I_1$ & irrational; in $\mathbb{Q}(\sqrt{5})$\cr
$\frac{11-5\sqrt{5}}{2}$ & $I_1$ & irrational; in $\mathbb{Q}(\sqrt{5})$\cr
\end{tabular}
\end{center}
\vspace{0.3cm}

As in the previous example, the double points of the fibre product $S_1(5) \times_{\mathbb{P}^{1}} Y'$ %%@
in the fibres above 0 and $\infty$ possess a projective small resolution. On the other hand, we %%@
will see later (section 6.2) that at least one of the two remaining nodes (in the irreducible %%@
singular fibres) allows no projective small resolution. Hence we resolve them by a usual blow-up %%@
each and denote the resulting mixed resolution by $\tilde{W}_2$. We will not consider the reduction of $\tilde{W}_2$ at 2 and 11, since $Y'$ is not well-defined there. Two further bad primes come from %%@
degeneration of singular fibres: In $\mathbb{F}_{5}$ we have again an $II\times II$ fibre instead of the %%@
two $I_1\times I_1$ fibres, requiring two further blow-ups. Meanwhile in characteristic 359 the %%@
congruence of $-\frac{125}{88}$ with one of the cusps $\frac{11\pm 5\sqrt{5}}{2}$ leads to an $III %%@
\times I_1$ fibre. At this point it is crucial to note that a small resolution would even resolve %%@
the degenerated singularity while our big resolution leaves one ordinary double point on the %%@
exceptional divisor such that we need one additional blow-up to resolve the degenerated %%@
singularity. Hence, the mixed resolution $\tilde{W}_2$ has bad reduction at 2, 5, 11 and 359.

\subsection{The fibre product $\tilde{W}_3$}

A third fibre product with one factor $S_1(5)$ and splitting third cohomology group can be derived %%@
from the elliptic surface $Y'$ with Weierstrass equation
\[
y^2=x(x-1)(x+(t^2-11t-1))
\]
since this gives for $t\in\mathbb{P}^{1}$ the following singular fibres:

\vspace{0.3cm}
\begin{center}
\begin{tabular}{c|c|c}
cusp & type of fibre & coordinates of nodes  \cr
\hline
$\infty$ & $I_4$ & rational\cr
$0$ & $I_2$ & rational\cr
$\frac{11+5\sqrt{5}}{2}$ & $I_2$ & irrational; in $\mathbb{Q}(\sqrt{-1})$\cr
$\frac{11-5\sqrt{5}}{2}$ & $I_2$ & irrational; in $\mathbb{Q}(\sqrt{-1})$\cr
$11$ & $I_2$ & rational
\end{tabular}
\end{center}
\vspace{0.3cm}

A mixed resolution $\tilde{W}_3$ of the fibre product $S_1(5)\times_{\mathbb{P}^{1}}Y'$ is obtained by %%@
blowing up the nodes of the two $I_1\times I_2$ fibres and some components of the two singular %%@
fibres over 0 and $\infty$. Concerning the reduction of $\tilde{W}_3$ we clearly have bad %%@
reduction at 2, since the Weierstrass equation of $Y'$ does not give elliptic fibres in this %%@
characteristic. Since $Y'$ has smooth reduction everywhere (as has $S_1(5)$), further bad primes %%@
can only come from degenerations of singular fibres. It is immediate that in characteristic 11 the %%@
elliptic surface $Y'$ has an $I_4$ fibre above 0, implying ten additional double points for the %%@
fibre product. On the contrary, although a similar thing happens in characteristic 5, this does %%@
not cause bad reduction, since the two $I_1 \times I_2$-fibres degenerate into one of type $II %%@
\times I_4$ such that the degenerated singularities are indeed resolved by the original four %%@
blow-ups.
Since reduction modulo other primes does not effect the fibrations, we conclude that $\tilde{W}_3$ %%@
has bad reduction only at 2 and 11.

Having examined the fibre products so far, we are now able to concentrate on the question of their %%@
modularity or, to be precise, of the modularity of the two-dimensional piece $U$ which we will %%@
isolate in their middle cohomology group in the next section.

\section{The Galois representations}

By the modularity of $U$ we mean that we can associate to $U$ a newform $f$ of weight 4 and level %%@
composed of $\tilde{W}$'s primes of bad reduction such that the good Euler factors of the %%@
$L$-series of $U$ coincide with those of the Mellin transform of $f$. Before turning our attention %%@
to a method for proving such an identity, let us first describe how to compute the Euler factors %%@
of the $L$-series of $U$.

The Galois representation $\varrho$ attached to $\cohom[\text{\'et}]{3}{\tilde{W}}{\mathbb{Q}_{\ell}}$ is %%@
known to be unramified outside the set $S$ consisting of $\ell$ and the primes of bad reduction of %%@
$\tilde{W}$. Hence for any good prime $p$, we can define the local Euler factor of %%@
$\cohom[\text{\'et}]{3}{\tilde{W}}{\mathbb{Q}_{\ell}}$ at $p$ as $\text{det}(1-\text{Frob}_p~p^{-s})^{-1}$. %%@
By the short exact sequence of (2) this product splits into factors corresponding to $U$ and $V$. %%@
This enables us to compute the $L$-series of $U$, with Euler factors at $p$ defined by
\[
\text{det}(1-\text{Frob}_p~p^{-s})\mid_U^{-1}=
(1-\text{trace Frob}_p\mid_Up^{-s}+\text{det Frob}_p\mid_Up^{-2s})^{-1},
\]
from the action of Frobenius on the complete cohomology group and on $V$.

Let at first $E\times\mathbb{P}^{1}$ be a component of a singular fibre contributing to $V$ (i.e. not %%@
meeting the 0-section). By the K\"unneth formula we have %%@
$\cohom[\text{\'et}]{3}{E\times\mathbb{P}^{1}}{\mathbb{Q}_{\ell}}=\cohom[\text{\'et}]{1}{E}{\mathbb{Q}_{\ell}}
\otimes\cohom[\text{\'et}]{2}{\mathbb{P}^{1}}{\mathbb{Q}_{\ell}}$, such that 
\[
\text{Frob}_p\mid_{\cohom[\text{\'et}]{3}{E\times\mathbb{P}^{1}}{\mathbb{Q}_{\ell}}}=p~
\text{Frob}_p\mid_{\cohom[\text{\'et}]{1}{E}{\mathbb{Q}_{\ell}}}.
\]

Since the elliptic curves $E$ will all be defined over $\mathbb{Q}$ in our examples, the results of Wiles %%@
et al. \cite{W} predict their modularity. Thus we find associated newforms $g_E$ of weight 2 and %%@
level the conductor of $E$. Hence, we have 
\[
L(V,s)=\prod_E L(g_E,s-1).
\]
On the other hand, the determinant of the Galois representation $\varrho$ of %%@
$\cohom[\text{\'et}]{3}{\tilde{W}}{\mathbb{Q}_{\ell}}$ is known to be the $(3\frac{h^3(\tilde{W})}{2})$-th %%@
power of the $\ell$-adic cyclotomic character, giving already $\text{det Frob}_p\mid_U=p^3.$ %%@
Meanwhile the trace of Frobenius on $\cohom[\text{\'et}]{3}{\tilde{W}}{\mathbb{Q}_{\ell}}$ can be computed %%@
by the Lefschetz fixed point formula, reading for a good prime $p$
\[
\#\tilde{W}(\mathbb{F}_{p})=\sum_{i=0}^{6} (-1)^i ~
\underbrace{\text{trace Frob}_p\mid_{\cohom[\text{\'et}]{i}{\tilde{W}}{\mathbb{Q}_{\ell}}}}
_{\text{tr}_i(p)}.
\]
Now it turns out that all summands but the trace tr$_3(p)$ are easily calculated. For this %%@
purpose, we use a computer to count the number of points over finite fields.  Looking at the right %%@
hand-term Poincar\'e duality implies that $\#\tilde{W}(\mathbb{F}_{p})$ equals %%@
$1+(1+p)\text{tr}_2(p)+p^3-\text{tr}_3(p)$ by the simple structure of the Hodge diamond. %%@
Furthermore, we conclude from section 2 that $\cohom[\text{\'et}]{2}{\tilde{W}}{\mathbb{Q}_{\ell}}$ is %%@
spanned by certain divisors (those generating Pic$\hat{W}$ and the exceptional quadrics). %%@
Since the operation of Frobenius on the divisor classes is easily derived from their defining %%@
equations, this allows us to compute the remaining trace tr$_3(p)$ for every good prime $p$ from %%@
the Lefschetz fixed point formula, eventually giving the local Euler factors of $L(U,s)$ at $p$ by %%@
virtue of
\[
\text{trace Frob}_p\mid_U=\text{tr}_3(p)-p\sum_E\text{trace %%@
Frob}_p\mid_{\cohom[\text{\'et}]{1}{E}{\mathbb{Q}_{\ell}}}.
\]
We will list some of the traces on $U$ as well as the number of points and the newforms associated %%@
to the elliptic curves $E$ in the last section. Then we look for a normalized newform %%@
$f=\sum_{n\in\mathbb{N}}a_nq^n$ of weight 4 in the tables of W. Stein \cite{St} with integer %%@
Fourier coefficients $a_p$ equaling the computed traces. This newform will indeed be appropriate %%@
since its Mellin transform $L(f,s)$ of $f$ has a product expansion with Euler factors %%@
$(1-a_pp^{-s}+p^{3-2s})^{-1}$ by the theory of Atkin-Lehner. The next step consists of applying a %%@
theorem of Deligne\cite{D} which associates to $f$ two-dimensional $\ell$-adic representations %%@
$\varrho_{f,\ell}$. Each is unramified outside the chosen $\ell$ and the prime divisors of the %%@
level $N$ of $f$ with determinant the third power of the $\ell$-adic cyclotomic character and %%@
trace $\varrho_{f,\ell}($Frob$_p)=a_p$ for all $p\nmid\ell N.$ This means that we can guarantee %%@
the equality of all good Euler factors of $L(U,s)$ and $L(f,s)$ by proving the isomorphicity of %%@
the semi-simplifications of the $\ell$-adic Galois representations $U$ and $\varrho_{f,\ell}$ for %%@
some $\ell$. Indeed this will be achieved by checking the explicit equality of only a small number %%@
of Euler factors. However, as the representations will not have even trace, this cannot be derived %%@
from the powerful theorem of R. Livn\'e \cite[Thm. 4.3]{L}.
Instead, we recapitulate an idea of J.-P. Serre, which to my knowledge has been applied only once %%@
before in this context by C. Schoen in \cite{Sc1}.

\section{Serre's construction}

In his \emph{R\'esum\'e des cours de 1984-1985}\cite{Se1}, J.-P. Serre explains a criterion to %%@
compare two-dimensional $2$-adic Galois representations, associated to elliptic curves. Here, we %%@
apply these ideas which are based on Faltings' work, beginning with a general construction %%@
which specializes to the above mentioned case of elliptic curves as well as to ours.

We start with two $\ell$-adic Galois representations
\[
\varrho_i: ~\text{Gal}(\bar{\mathbb{Q}}/\mathbb{Q}) \rightarrow \text{GL}_n(\mathbb{Q}_{\ell})~~~~(i=1,2)
\] 
which both are unramified outside a finite set of primes $S$. After specifying a stable lattice we can as well assume them to take values in GL$_n(\Z_{\ell})$ (cf. \cite[1.1]{Se2}]. Then we make the assumption that $\varrho_1$ and $\varrho_2$ fulfill the following conditions:

\begin{enumerate}
\item They have the same determinant.
\item The mod $\ell$ reductions $\bar{\varrho}_1, \bar{\varrho}_2$ are absolutely irreducible and %%@
isomorphic.
\end{enumerate}

This obviously gives the congruence tr $\varrho_1\equiv$ tr $\varrho_2$ mod $\ell$. Now assume %%@
further that the traces of $\varrho_1$ and $\varrho_2$ are not identical, i.e. there is some prime %%@
$p\not\in S$ such that tr $\varrho_1($Frob$_p)\neq$ tr $\varrho_2($Frob$_p)$. Then, choose the %%@
maximal $\alpha\in\N$ such that\linebreak
tr $\varrho_1\equiv$ tr $\varrho_2$ mod $\ell^{\alpha}$ (or, equivalently, $\varrho_1$ and $\varrho_2$ are isomorphic mod $\ell^{\alpha}$). Hereby we obtain a non-constant map
\begin{eqnarray*}
\tau: ~\text{Gal}(\bar{\mathbb{Q}}/\mathbb{Q}) & \rightarrow & ~~~~~~~~~~~\mathbb{F}_{\ell}\\
\sigma~~~~~~~ & \mapsto & \frac{\text{tr}~ \varrho_1(\sigma)-\text{tr}~ %%@
\varrho_2(\sigma)}{\ell^{\alpha}} ~~\text{mod}~ \ell
\end{eqnarray*}
which maps the inertia groups $I_p$ to 0 for $p\not\in S$. Our next aim is to construct a suitable %%@
factorization $\tau=\widetilde{\tau}\circ\widetilde{\varrho}$ and then investigate the map %%@
$\widetilde{\tau}$.

For this purpose -after replacing $\varrho_1$ by a conjugate if necessary- we can assume that $\varrho_1\equiv \varrho_2$ mod $\ell^{\alpha}$. Hence for %%@
every $\sigma\in\text{Gal}(\bar{\mathbb{Q}}/\mathbb{Q})$ there is a matrix $\mu(\sigma)\in$ M$_n(\Z_{\ell})$ with
\[
\varrho_1(\sigma)=(1+\ell^{\alpha}~\mu(\sigma))\varrho_2(\sigma).
\]
As this relation allows us to express $\tau$ as $\tau(\sigma)=$ tr $\mu(\sigma)\varrho_2(\sigma),$ %%@
we are lead to factor $\tau$ through the product M$_n(\Z_{\ell})\times$ GL$_n(\Z_{\ell})$. At this %%@
point we notice that due to the definition of $\tau$ mod $\ell$ we can restrict ourselves to the %%@
product of the mod $\ell$ reductions $\bar{\mu}$ and $\bar{\varrho}_2$ (or, say, $\bar{\varrho}$ %%@
~due to condition 2.). Then, the map $\widetilde{\varrho}=\bar{\mu} \times \bar{\varrho}$ is made %%@
into a group homomorphism by giving the target group M$_n(\mathbb{F}_{\ell}) \times$ GL$_n(\mathbb{F}_{\ell})$ the %%@
structure of a semi-direct product M$_n(\mathbb{F}_{\ell}) \rtimes$ GL$_n(\mathbb{F}_{\ell})$ with operation 
\[
(A,C)\cdot(B,D)=(A+CBC^{-1},CD).
\] 
By construction, $\widetilde{\varrho}$ is unramified outside $S$. Furthermore, condition 1. %%@
implies that det$(1+\ell^{\alpha}~\mu)=1$. Expanding this equation in terms of powers of $\ell$, %%@
this equation becomes $1=1+\ell^{\alpha}$ tr $\mu +\ell^{2\alpha}(...)$. Hence, tr $\mu\equiv 0$ %%@
mod $\ell$ and 
\[
\widetilde{\varrho}:~ \text{Gal}(\bar{\mathbb{Q}}/\mathbb{Q}) \rightarrow \widetilde{G}
\]
where $\widetilde{G}$ is the subgroup of the semi-direct product M$_n(\mathbb{F}_{\ell}) \rtimes$ %%@
GL$_n(\mathbb{F}_{\ell})$ consisting of all elements whose first component have trace 0 mod $\ell$. In the %%@
following we apply this method, measuring the discrepancy between $\varrho_1$ and $\varrho_2$ by %%@
mapping the absolute Galois group into a finite group via $\widetilde{\varrho}$, to our special %%@
situation:

Let $n=2$ and $\ell=2.$ We construct the map $\widetilde{\varrho}$ in the case where both, %%@
$\varrho_1$ and $\varrho_2$, do not have even trace, but the same determinant $\chi_2^3$.
Then condition 2. is easily checked by computing the Galois extensions $K_i/\mathbb{Q}$ cut out by the %%@
kernels of the mod 2 reductions $\bar{\varrho}_i$ for $i=1,2$. Indeed, the absolute irreducibility %%@
of $\bar{\varrho}_i$ is equivalent to $K_i/\mathbb{Q}$ having Galois group $S_3$ in our situation, since %%@
the traces being not even imply the Galois extensions $K_i/\mathbb{Q}$ to have Galois group $S_3$ or $C_3$ %%@
(the image of $\bar{\varrho}_i$ in GL$_2(\mathbb{F}_{2})$). Subsequently, $\bar{\varrho}_1$ and %%@
$\bar{\varrho}_2$ are furthermore isomorphic if $K_1=K_2$. Explicitly, we proceed in the following %%@
way: As the $K_i$ are unramified outside $S$, there is only a finite number of isomorphism classes %%@
of them which can be computed by class field theory. We find the possible number fields in the %%@
tables of J. Jones \cite{J}. Then the number field $K_i$ cut out by $\bar{\varrho}_i$ can be %%@
determined by the fact that Frob$_p$ has order 3 in Gal$(K_i/\mathbb{Q})$ if and only if its trace is odd. %%@
Thereby, the simultaneous oddness of the traces tr $\varrho_i(\text{Frob}_p)$ for some suitable %%@
primes $p$ guarantees the isomorphicity of the mod 2 reductions $\bar{\varrho}_1$ and %%@
$\bar{\varrho}_2$. 

In this setting, condition 1. and 2. are fulfilled, such that we can construct the map %%@
$\widetilde{\varrho}$ as sketched above under the assumption that the $\varrho_i$ be not
isomorphic. Next we will investigate the structure of the target group %%@
$\widetilde{G}$ of $\widetilde{\varrho}$. By inspection $\widetilde{G}$ has at most 48 elements. %%@
Indeed, our next aim is to prove that $\widetilde{G}\cong S_4 \times C_2$. In order to see this, %%@
consider the bijection 
\begin{eqnarray*}
j:~~~~~~~~~~ \widetilde{G}~~~~~~~~~ & \rightarrow & ~~~(\mathbb{F}_{2}^{~2} \times \text{GL}_2(\mathbb{F}_{2})) %%@
\times \mathbb{F}_{2}\\
\left( \begin{pmatrix} a_1 & a_2\\ a_3 & a_1\end{pmatrix}, C\right) & \mapsto & %%@
\left(\left(\begin{pmatrix}a_2\\ a_3\end{pmatrix}, C\right), a_1+a_2+a_3\right). 
\end{eqnarray*}
This becomes a group isomorphism relative to the group structure of $\widetilde{G}$ if we define %%@
the group operation on the target group of $j$ as follows. For the first two factors we define %%@
multiplication as
\[
\left(\begin{pmatrix}a_2\\ a_3\end{pmatrix}, C\right)\cdot\left(\begin{pmatrix}b_2\\ %%@
b_3\end{pmatrix}, D\right) = \left(\begin{pmatrix}a_2\\ a_3\end{pmatrix}+C\begin{pmatrix}b_2\\ %%@
b_3\end{pmatrix}, CD\right)
\]
while in the third factor we take the addition on $\mathbb{F}_{2}$. By virtue of $j$, we therefore identify %%@
$\widetilde{G}$ with a product of a group $\widetilde{G}'$ with $\mathbb{F}_{2}$. By inspection %%@
$\widetilde{G}'$ contains 24 elements including eight order 3 and six of order 4. Hence, it has to %%@
coincide with $S_4$, giving the claimed isomorphism $\widetilde{G}\cong S_4 \times C_2$.

At this point, recall that the assumption that $\varrho_1$ and $\varrho_2$ be not isomorphic %%@
implies the map $\tau$ to be non-constant. In order to establish a %%@
contradiction, we consider the subsequent map 
\begin{eqnarray*}
\widetilde{\tau}: ~~~\widetilde{G}~~~ & \rightarrow &~~~ \mathbb{F}_{2}\cr
(A,C) & \mapsto & \text{tr}~ A \cdot C~ \text{mod}~ 2
\end{eqnarray*}
which gives $\tau=\widetilde{\tau}\circ\widetilde{\varrho}$. Direct calculations show that %%@
$\widetilde{\tau}$ does not vanish exactly at the elements of $\widetilde{G}$ of order greater %%@
than 3. As a consequence, the assumption forces the image of $\widetilde{\varrho}$ to contain an %%@
element of order 4 or 6. Then, by the structure of the target group $\widetilde{G}$ the image of %%@
$\widetilde{\varrho}$ is already seen to include $S_3\times C_2$ or $S_4$ (since im %%@
$\bar{\varrho}=S_3$ is clearly contained in im $\widetilde{\varrho}$). By %%@
Galois theory, the image of $\widetilde{\varrho}$ corresponds to a Galois extension $L/\mathbb{Q}$. By the %%@
construction of $\widetilde{\varrho}$ the number field $L$ is unramified outside $S$ with %%@
intermediate field, say, $K=K_1=K_2$ (the extension cut out by ker$\bar{\varrho}$). Here, the main %%@
point is that Galois extensions with small Galois group and ramification locus have been %%@
classified to some extent (cf. \cite{J}). Indeed, in our situation we can restrict $L/\mathbb{Q}$ to have %%@
Galois group $S_3\times C_2$ or $S_4$ by the above considerations. Then, having found every %%@
possible such extension thanks to the calculations of J. Jones, we will rule out each one of these %%@
by an easy method: Find Frob$_p$ with maximal order in Gal$(L/\mathbb{Q})$; this implies %%@
$\widetilde{\tau}($Frob$_p)=1$. Then check the explicit equality\linebreak
tr $\varrho_1($Frob$_p)=$tr $\varrho_2($Frob$_p)$, giving the desired contradiction.

To sum it up the proof of the isomorphicity of $\varrho_1$ and %%@
$\varrho_2$ amounts to checking the explicit equalities of the traces of Frob$_p$ for a small %%@
finite number of primes. We will execute this procedure for our examples in the next section %%@
exactly as sketched here. Note, however, that we may apply these ideas as well without determining %%@
the intermediate field $K$ explicitly. The only difference is that this requires to carry out the %%@
whole procedure for any cubic extension with Galois group $S_3$ which is unramified outside $S$. %%@
Eventually this thereby leads to a quite general criterion for the isomorphicity of %%@
(semi-simple) two-dimensional $2$-adic Galois representations.

The notably easiest example of this kind is implicitly included in a letter from Serre to S. Bloch %%@
(cf. \cite{Sc1}) and the paper of Livn\'e \cite{L}. It says that a semi-simple Galois representation with %%@
ramification locus $\{2,5\}$ is uniquely determined by its determinant and the traces of Frobenius %%@
at $p\in T=\{3,7,11,13,17,29,31\}$. Here we will treat another example of this kind which will be %%@
used for the final example of this article:

Let $S=\{2,3\}$ and $\varrho_1,\varrho_2$ be two-dimensional 2-adic Galois %%@
representations with the same determinant and ramification locus $S$. Start with the tables of J. %%@
Jones \cite{J} to find the nine number fields with Galois group $S_3$ or $C_3$ and ramification %%@
locus $S$. By inspection of the Galois groups, the traces of $\bar{\varrho}_i($Frob$_p)$ for $p\in %%@
T=\{5,7,11,13,31\}$ determine uniquely whether and which one is cut out by ker $\bar{\varrho}_i$. %%@
Note that the $C_3$-extension is the only one with odd traces at both, 11 and 13. Adding 17, 19, %%@
and 37 to $T$ guarantees that for any possible $S_3\times C_2$ or $S_4$ extension, there is %%@
at least one $p\in T$ with Frob$_p$ of order 6 or 4, respectively. Adding furthermore 23 to $T$, %%@
we finally obtain a set which is sufficient in Livn\'e's sense, as the compositum of all quadratic %%@
extensions which are unramified outside $S$ is directly seen to be
$Q(\zeta_{24})$ with $\zeta_{24}$ a primitive 24-th root of unity and $\{\text{Frob}_p: %%@
p=5,7,11,13,17,19,23\}=\text{Gal}(\mathbb{Q}(\zeta_{24})/\mathbb{Q})-\{1\}.$ Hence we have the following

\textbf{Proposition 1.}\newline
Let $\varrho_1, \varrho_2$ be two-dimensional 2-adic Galois representations with the same %%@
determinant and even trace at Frob$_{11}$ or Frob$_{13}$, which are unramified outside $\{2,3\}.$ %%@
Then they have isomorphic semi-simplifications if and only if for every \linebreak%%@
$p\in\{5,7,11,13,17,19,23,31,37\}$
\[
\text{trace}~ \varrho_1(\text{Frob}_p)=\text{trace}~ \varrho_2(\text{Frob}_p).
\]

We will use this proposition to prove the modularity of another Calabi-Yau threefold to be %%@
introduced in the very last section, but for the moment we return to our three fibre products.

\section{The modularity of the fibre products}

We are going to apply the ideas of the last section to the 2-adic Galois representations $U$ and %%@
$\varrho_{f,2}$ described in section 4 where $U$ lives in the third cohomology group of one of the %%@
fibre products $\tilde{W}_i$. For this purpose, we collect the remaining properties of the fibre %%@
products necessary to compute the traces of Frobenius on $U$. Then we compute a short list of %%@
traces. By the calculations of J. Jones \cite{J} these will indeed turn out to be sufficient in %%@
the sense of Serre's approach except for one case. Hence we look in the tables of W. Stein %%@
\cite{St} for a newform of weight 4 with exactly these traces as coefficients and express the %%@
$L$-series and the zeta function of the fibre products by means of this newform. Throughout we %%@
employ the notation of the previous sections.

\subsection{The variety $\tilde{W}_1$}

The threefold $\tilde{W}_1$ has $h^{1,1}=37$ and $h^{1,2}=8$. The second equality, which follows %%@
from equation (1), can be derived from $\tilde{W}_1$ having two fibres of type $I_0\times I_5$ %%@
above 0 and 11. Meanwhile the first equality follows immediately from the considerations about the %%@
Picard group. More precisely, the 37 divisors spanning %%@
$\cohom[\text{\'et}]{2}{\tilde{W}_1}{\mathbb{Q}_{\ell}}$ are directly seen to be all defined over $\mathbb{Q}$ %%@
except for the two exceptional quadrics. Since these are exchanged by Frob$_p$ if and only if 5 is %%@
not a square modulo $p$, we have tr$_2(p)=(36+(\frac{5}{p}))p.$ On the other hand, it is clear %%@
from the previous observations that all of the 16-dimensional quotient $V$ of %%@
$\cohom[\text{\'et}]{3}{\tilde{W}_1}{\mathbb{Q}_{\ell}}$ is made up by eight copies of the twisted Galois %%@
representation $\cohom[\text{\'et}]{1}{E}{\mathbb{Q}_{\ell}}(1)$ where $E$ is the fibre of $S_1(5)$ above %%@
11. This elliptic curve has discriminant $-11^5$ and conductor $11$ and is thus associated to the %%@
unique newform of weight 2 and level 11, 
\[
g_{11}=q - 2q^2 - q^3 + 2q^4 + q^5 + 2q^6 - 2q^7 - 2q^9 - 2q^{10} + q^{11} - 2q^{12} + %%@
4q^{13}\hdots
\]
This enables us to compute the traces of Frob$_p$ on $U$ from the number of points of %%@
$\tilde{W}_1$ over $\mathbb{F}_{p}$ as indicated in the following table:

{\scriptsize
\begin{center}
\begin{tabular}{c|c|c}
$p$ &  $\sharp\tilde{W}_1(\mathbb{F}_{p})$ & tr$_U(p)$\cr
\hline
3 & 475 & -3\cr
7 & 2425 & -9\cr
13 & 8150 & 2\cr
17 & 15875 & 21\cr
23 & 31650 & 22\cr
41 & 135738 & -478\cr
43 & 147800 & -8
\end{tabular}
\end{center}}

Since $\tilde{W}_1$ has bad reduction only at 5 and 11, $U$ is unramified outside $S=\{2,5,11\}.$
By the tables of J. Jones \cite{J} the oddness of the traces for 3, 7 and 17 determines the Galois %%@
extension cut out by the mod 2 reduction of $U$ as $\mathbb{Q}(x^3+2x-8)/\mathbb{Q}$ with Galois group $S_3$ and %%@
intermediate field $\mathbb{Q}(x^2+110)$. The other quadratic extensions of $\mathbb{Q}$ which are unramified %%@
outside $S$ each have generating polynomial irreducible mod 3, 7 or 17 (such that Frob$_3$, %%@
Frob$_7$ or Frob$_{17}$ has order 2 in the corresponding Galois group). Hence, these three primes %%@
are sufficient to prevent an element of order 6 in the image of the representation $\widetilde{U}$ %%@
(here $\widetilde{U}$ is constructed like $\widetilde{\varrho}$ in the last section).

In order to find sufficient primes to exclude the elements of order 4 as well, J. Jones was so %%@
kind as to run his programs for the Galois extensions $L/\mathbb{Q}$ with Galois group $S_4$ and %%@
intermediate field $\mathbb{Q}(x^3+2x-8)$. By his calculations these have generating polynomials
{\small 
\begin{center}
\begin{tabular}{ccc}
$x^4 - 220x + 165$&$x^4 - 2x^3 - 26x^2 - 28x - 24$&$x^4 - 10x^2 - 40x + 10$ \cr
$x^4 - 2x^3 - x^2 - 8x - 4$&$x^4 - 2x^3 - 15x^2 - 28x - 24$&$x^4 - 22x^2 - 176x - 110$ \cr
$x^4 - 440x - 3410$&$x^4 - 20x^2 - 40x - 10$&$x^4 + 20x^2 - 80x + 60$ \cr
$x^4 + 44x^2 - 352x + 132$&$x^4 - 12x^2 - 16x - 20$&$x^4 - 8x - 2$ \cr
$x^4 - 4x^2 - 4x + 9$&$x^4 - 176x + 418$&$x^4 - 1760x - 440.$
\end{tabular}
\end{center}}
Each of these polynomials is irreducible modulo 13, 23, 41 or 43. As a consequence, for each %%@
Galois extension we find Frob$_p$ for one of these four primes with order 4 in the corresponding %%@
Galois group. Together with 3,7 and 17 this gives the sufficient set of seven primes for which we %%@
listed the traces in the table above. Hence, if we can find a newform $f$ of weight 4 with Fourier %%@
coefficients these traces, Serre's method will enable us to conclude the representations $U$ and %%@
$\varrho_{f,2}$ to be isomorphic and thereby the modularity of $U$. Indeed, %%@
the tables of W. Stein \cite{St} give the unique newform of level 55 with rational coefficients, %%@
$f_{55}=q + q^2 - 3q^3 - 7q^4 - 5q^5 - 3q^6 - 9q^7 - 15q^8 - 18q^9 - 5q^{10} + 11q^{11} + 21q^{12} %%@
+ 2q^{13}\hdots .$ Hereby, with $\chi_5$ the Dirichlet character cutting out $\mathbb{Q}(\sqrt{5})$ and %%@
$\circeq$ denoting equality up to a finite number of Euler factors, we have proven the

\textbf{Theorem 2.}\newline
The variety $\tilde{W}_1$ is modular. We can isolate a two-dimensional piece $U$ in its third %%@
cohomology group whose good Euler factors coincide with those of the newform $f_{55}$ of weight 4 %%@
and level 55. $L$-series and $\zeta$-function are given as follows:
\[
L(U,s)\circeq L(f_{55},s),~~~~~~~~ L(\tilde{W}_1,s)\circeq L(g_{11},s-1)^8 L(f_{55},s)~~~~~~~~~~ \text{and} %%@
\]\[
\zeta(\tilde{W}_1,s)\circeq\frac{\zeta(s)\zeta(s-1)^{36}
\zeta(\chi_5,s-1)\zeta(s-2)^{36}\zeta(\chi_5,s-2)\zeta(s-3)}{L(\tilde{W}_1,s)}.
\]

\subsection{The variety $\tilde{W}_2$}

This variety has $h^{1,1}=45$ and $h^{1,2}=1$, again by equation (1) and the observations of the %%@
following paragraph. Furthermore, the 45 divisors spanning the Picard group of  $\tilde{W}_2$ %%@
behave very similarly to the above example (i.e. only the two exceptional divisors are not eigen %%@
elements for all Frobenii). The only divisors for which this is not trivial are those generating %%@
the rank 2 Picard group of the generic fibre of the elliptic surface $Y'$. However, by the %%@
description as a cubic pencil, we get $Y'$ as $\mathbb{P}^{2}$ blown up in the nine base points of this %%@
pencil. Since these are all rational, we have $\# Y'(\mathbb{F}_{p})=1+10p+p^2$, such that Frob$_p$ has to %%@
operate as multiplication by $p$ ~on $\cohom[\text{\'et}]{2}{Y'}{\mathbb{Q}_{\ell}}$ (and in particular on %%@
Pic $Y'_{\eta}$) by the Lefschetz fixed point formula. We know the two-dimensional representation %%@
$V$ in $\cohom[\text{\'et}]{3}{\tilde{W}_2}{\mathbb{Q}_{\ell}}$ to be generated by the fibre $E$ of %%@
$S_1(5)$ above $-\frac{125}{88}$. This elliptic curve has discriminant %%@
$-5^{15}359^22^{-21}11^{-7}$ and conductor $2\cdot 5\cdot 11\cdot 359$. The coefficients of the %%@
corresponding weight 2 newform $f_{39490}=q + q^2 - q^3 + q^4 + q^5 - q^6 + 3q^7 + q^8 - 2q^9 + %%@
q^{10} + q^{11} - q^{12} + 4q^{13} \hdots$ eventually enable us to compute the following data:

{\scriptsize
\begin{center}
\begin{tabular}{c|c|c}
$p$ &  $\sharp\tilde{W}_2(\mathbb{F}_{p})$ & tr$_U(p)$\cr
\hline
3 & 550 & -3\cr
7 & 2740 & -9\cr
13 & 9970 & 2\cr
17 & 18000 & 21\cr
23 & 35790 & 22\cr
41 & 147218 & -478\cr
43 & 160700 & -8
\end{tabular}
\end{center}}

One could of course compute more traces, but we will content ourself with these few, as they all %%@
coincide with those calculated for $\tilde{W}_1$. However, the threefold $\tilde{W}_2$ has bad %%@
reduction at 359, such that the data above will a priori be not sufficient (unlike it was if %%@
$\tilde{W}_2$ had bad reduction only at 2, 5 and 11). Unfortunately there does not seem to be a %%@
method at hand to prove isomorphicity of Galois representations in this case, since the bad prime %%@
359 appears to be much too large to compute all cubic and quartic Galois extensions unramified %%@
outside $\{2,5,11,359\}$ (at least at the time being). Although it seems reasonable that the ramification at 359 takes only place on $V$, not on $U$, it is only speculation that the regular three form should still exist on a suitable model of the special fibre. Hence, we can only stress the numerical %%@
evidence and formulate the following

\textbf{Conjecture 3.}\newline
The variety $\tilde{W}_2$ is modular. We can isolate a two-dimensional piece $U$ in its third %%@
cohomology group whose good Euler factors coincide with those of the newform $f_{55}$ of weight 4 %%@
and level 55. $L$-series and $\zeta$-function are given as follows:
\[
L(U,s)\circeq L(f_{55},s),~~~~~~~~ L(\tilde{W}_2,s)\circeq L(g_{39490},s-1) L(f_{55},s)~~~~~~~~~~ \text{and} %%@
\]\[
\zeta(\tilde{W}_2,s)\circeq\frac
{\zeta(s)\zeta(s-1)^{44}\zeta(\chi_5,s-1)\zeta(s-2)^{44}\zeta(\chi_5,s-2)\zeta(s-3)}{L(\tilde{W}_2,s)}.
\]

One interesting point about this conjecture is that there are bad primes of $\tilde{W}_2$ which do %%@
not appear in the level of $f_{55}$. To my knowledge this is the first example where such a thing %%@
is conjectured to happen, so it might be worth spending a closer look. We saw in section 3.2 that %%@
a projective small resolution, did it exist, would have good reduction at 359. However, the %%@
conductor of $V=\cohom[\text{\'et}]{1}{E}{\mathbb{Q}_{\ell}}(1)$ contradicts the projectivity of any small %%@
resolution since it predicts the Galois representation of its third cohomology group, decomposing %%@
into the two-dimensional pieces $V$ and $U$, to be ramified at 359. Nevertheless, the different %%@
behaviour of small and big resolution seems to be reflected in the level of the newform which we %%@
conjecture to correspond to $U$. It should be stressed that this does not seem to be an %%@
exceptional case. Indeed, there are some more examples of nodal non-rigid Calabi-Yau threefolds %%@
with similar behaviour which also arise from fibre products of elliptic surfaces and which will %%@
hopefully be subject of a future work. But at the moment, let us now turn to our third fibre %%@
product $\tilde{W}_3$.

\subsection{The variety $\tilde{W}_3$}

By the same considerations as above $\tilde{W}_3$ has $h^{1,1}=39$ and $h^{1,2}=1$. But unlike %%@
before the computation of the action of Frobenius on the second cohomology group %%@
requires some more work in this case. This is due to the rank 2 Picard group of the generic fibre of $Y'$ together with the fact that we did not realize $Y'$ as $\mathbb{P}^{2}$ blown up in the nine base points of a cubic pencil. However, we can read off eight of the ten eigenvalues of Frob$_p$ on $\cohom[\text{\'et}]{2}{Y'}{\mathbb{Q}_{\ell}}$ to be seven times $p$ and once $\chi_5(p) p$ directly from the components. Thus the eigenvalues on Pic $Y'_{\eta}$ have to be $p$ and $\pm p$ by the Weil conjectures, with the sign of the second given by a quadratic character which is unramified outside 2 and 5. Since there are only seven such characters, we can specify the actual one by a few explicit calculations. Indeed the sign of the last eigenvalue turns out to be $-1$ for $p=3,7,13$, and the only character of the above which satisfies this property is given by $\chi_5$.

As before, the action of Frobenius on the remaining part of the second cohomology group is easily understood: The 0-section and the components of the fibres %%@
above $0, 11$ and $\infty$ are all rational. The components of the two $I_1\times I_2$ fibres %%@
above $\frac{11\pm 5\sqrt{5}}{2}$ are exchanged by Frob$_p$ if and only if 5 is a quadratic %%@
non-residue in $\mathbb{F}_{p}$. Similarly, the four exceptional quadrics are permuted by the Klein %%@
four-group (i.e. they are only fixed by Frob$_p$ if both, 5 and -1, are squares modulo $p$). %%@ 
Finally we know the fibre of $S_1(5)$ above 11, whose %%@
first cohomology group generates the two-dimensional quotient $V$ of %%@
$\cohom[\text{\'et}]{3}{\tilde{W}_3}{\mathbb{Q}_{\ell}}$ after a Tate twist, to be associated to $g_{11}$ %%@
from section 6.1, so we can compute

{\scriptsize
\begin{center}
\begin{tabular}{c|c|c}
$p$ &  $\sharp\tilde{W}_3(\mathbb{F}_{p})$ & tr$_U(p)$\cr
\hline
3 & 410 & -7\cr
7 & 2080 & 14\cr
13 & 7860 & -72\cr
23 & 29410 & -107\cr
31 & 64178 & 117
\end{tabular}
\end{center}}

It will turn out that this small set of data is sufficient to prove the modularity of $U$. Note however, that a priori we only 
know the Galois representation $U$ to be unramified outside $S=\{2,5,11\}$, since the conjugation of $\mathbb{Q}(\sqrt{5})/\mathbb{Q}$ acts non-trivially on $\cohom[\text{\'et}]{2}{\tilde{W}_3}{\mathbb{Q}_{\ell}}$ by the above inspections. By the tables of \cite{J} the Galois extension %%@
$\mathbb{Q}(x^3-x^2+x+1)/\mathbb{Q}$, cut out by the kernel of the mod 2 reduction of $U$, is uniquely determined %%@
by the simultaneous oddness of tr$_U(3)$ and evenness of tr$_U(7)$. While it has Galois group $S_3$ and intermediate field $\mathbb{Q}(x^2-x+3)$, %%@
all other quadratic extensions which are unramified outside $S$ have Frob$_p$ of order 2 in their %%@
Galois group for $p=3,23$ or $31$ (and $x^3-x^2+x+1$ is irreducible modulo $3,23$ and $31$). This %%@
settles the case of $S_3 \times C_2$ extensions. On the other hand, there are only three quartic Galois %%@
extensions of $\mathbb{Q}$ with Galois group $S_4$ and intermediate field $\mathbb{Q}(x^3-x^2+x+1)$, which are unramified outside $S$, by calculations of J. Jones. Indeed, these are furthermore unramified at 5, having discriminant -11, and can thus be found at \cite{J}. Since in each %%@
case Frob$_7$ or Frob$_{13}$ has order 4 in the Galois group, the primes 3, 7, 13, 23 and 31 are %%@
indeed sufficient for Serre's approach. Finding the traces as Fourier coefficients of a newform of %%@
weight 4 and level 22, $f_{22}=q - 2q^2 - 7q^3 + 4q^4 - 19q^5 + 14q^6 + 14q^7 - 8q^8 + 22q^9 + %%@
38q^{10} + 11q^{11} - 28q^{12} - 72q^{13} \hdots$, we thereby have proven the modularity of %%@
$U$.

\textbf{Theorem 4.}\newline
The variety $\tilde{W}_3$ is modular. We can isolate a two-dimensional piece $U$ in its third %%@
cohomology group whose good Euler factors coincide with those of the newform $f_{22}$ of weight 4 %%@
and level 22. The $L$-series and zeta-function are given as follows:
\[
L(U,s)\circeq L(f_{22},s),~~~~~~~~ L(\tilde{W}_3,s)\circeq L(g_{11},s-1) L(f_{22},s)~~~~~~~~~~\text{and}\]

\[
\zeta(\tilde{W}_3,s)\circeq  \frac{\begin{matrix}\zeta(s)\zeta(s-1)^{34}
\zeta(\chi_{5},s-1)^{3}\zeta(\chi_{-1},s-1)\zeta(\chi_{-5},s-1)~~~~~~~~\\
~~~~~~\zeta(s-2)^{34}\zeta(\chi_{5},s-2)^{3}\zeta(\chi_{-1},s-2)\zeta(\chi_{-5},s-2)\zeta(s-3)\end{matrix}}{L(\tilde{W}_3,s)}.
\]

\subsection{A Calabi-Yau threefold, associated to a newform of level 27}

As announced in the introduction, we would eventually like to present another example of non-rigid %%@
modular Calabi-Yau threefolds, whose modularity will be proven by virtue of Proposition 1. It %%@
arises from the Hesse-pencil 
\[
x^3+y^3+z^3+txyz=0
\] 
which also appears in Beauville's classification\cite{B}. This modular elliptic surface $S(3)$, %%@
associated to the principal congruence subgroup $\Gamma(3)$, has four singular fibres of type %%@
$I_3$ over the cusps $\infty, -3, -3\omega$ and $-3\omega^2$
where $\omega$ denotes a primitive third root of unity. Hence, modifying the projection onto %%@
$\mathbb{P}^{1}$ by the automorphism $\pi: t\mapsto 3-t$ allows us to form a fibre product $W=(S(3), %%@
pr)\times_{\mathbb{P}^{1}}(S(3),\pi\circ pr)$ with singularities above the three cusps $\infty, %%@
-3\omega$ and $-3\omega^2$. Since all components of the singular fibres are smooth surfaces, it %%@
possesses a projective small resolution $\hat{W}$. Similar to the other three examples we see that %%@
the right hand-side of equation (1) reduces to the last summand and conclude $h^{1,1}=31$ and %%@
$h^{1,2}=4$. Thus, all of the 8-dimensional quotient $V$ of %%@
$\cohom[\text{\'et}]{3}{\hat{W}}{\mathbb{Q}_{\ell}}$ is generated by four copies of %%@
$\cohom[\text{\'et}]{3}{E}{\mathbb{Q}_{\ell}}(1)$ where $E$ is the fibre of $S(3)$ above $6$ which appears %%@
in the two $I_3 \times I_0$ fibres above 0 and 6. This elliptic curve has discriminant %%@
$-2^{12}3^3$ and conductor $27$ with associated newform $g_{27}=q - 2q^4 - q^7 + 5q^{13} \hdots$. %%@
Finally we note that only the fibre above $\infty$ and the divisor classes of the 0-section and %%@
the generic fibre are defined over $\mathbb{Q}$. The 20 other divisors spanning %%@
$\cohom[\text{\'et}]{2}{\hat{W}}{\mathbb{Q}_{\ell}}$ are only defined over $\mathbb{Q}(\sqrt{-3})$ and are %%@
therefore permuted by Frob$_p$ if $-3$ is not a square in $\mathbb{F}_{p}$. As this gives %%@
tr$_2(p)=(21+10(\frac{-3}{p})p)$, we can thereby compute the trace of Frobenius on $U$ (the %%@
two-dimensional orthogonal complement of $V$ in $\cohom[\text{\'et}]{3}{\hat{W}}{\mathbb{Q}_{\ell}}$) as %%@
follows:

{\scriptsize
\begin{center}
\begin{tabular}{c|c|c}
$p$ &  $\sharp\hat{W}(\mathbb{F}_{p})$ & tr$_U(p)$\cr
\hline
5 & 471 & -15\cr
7 & 2133 & -25\cr
11 & 2769 & 15\cr
13 & 7560 & 20\cr
17 & 8352 & -72\cr
19 & 19170 & 2\cr
23 & 18354 & -114\cr
31 & 60939 & 101\cr
37 & 93042 & -430
\end{tabular}
\end{center}}

We find these traces as coefficients of a newform of weight 4 and level 27, $f_{27}=q - 3q^2 + q^4 %%@
- 15q^5 - 25q^7 + 21q^8 + 45q^{10} + 15q^{11} + 20q^{13} \hdots$. Since $\hat{W}$ has bad %%@
reduction exactly at 3, we can apply Proposition 1 to deduce the isomorphicity of the 2-adic Galois representations $U$ and $\varrho_{f_{27},2}$:

\textbf{Theorem 5.}\newline
The Calabi-Yau threefold $\hat{W}$ is modular. We can isolate a two-dimensional piece $U$ in its %%@
third cohomology group whose good Euler factors coincide with those of the newform $f_{27}$ of %%@
weight 4 and level 27. The $L$-series and $\zeta$-function are given as follows with $\chi_{-3}$ %%@
the Dirichlet character cutting out $\mathbb{Q}(\sqrt{-3})$:
\[
L(U,s)\circeq L(f_{27},s),~~~~~~~~ L(\hat{W},s)\circeq L(g_{27},s-1)^4 L(f_{27},s)~~~~~~~~~~\text{and}\]\[
\zeta(\hat{W},s)\circeq\frac{\zeta(s)\zeta(s-1)^{21}
\zeta(\chi_{-3},s-1)^{10}\zeta(s-2)^{21}\zeta(\chi_{-3},s-2)^{10}\zeta(s-3)}{L(\hat{W},s)}.
\]

Matthias Sch\"utt,\newline
Institut f\"ur Mathematik (C),\newline
Universit\"at Hannover, \newline
Welfengarten 1, \newline
30060 Hannover, \newline
Germany,\newline
{\tt schuett@math.uni-hannover.de}


{\small
\begin{thebibliography}{9999}
\bibitem[B]{B} Beauville, A.: \emph{Les familles stables de courbes elliptiques sur $\mathbb{P}^{1}$ %%@
admettant quatre fibres singulière}, C.R. Acad. Sci. Paris 294 (1982), 657-660.
\bibitem[D]{D} Deligne, P.: \emph{Formes modulaires et repr\'esentations $\ell$-adiques}, in: %%@
S\'em. Bourbaki 1968/69, no. 355 (Lect. Notes in Math. 179), Springer-Verlag (1971), 139-172.
\bibitem[DM] {DM} Dieulefait, L., Manoharmayum, J.: \emph{Modularity of rigid Calabi-Yau %%@
threefolds over $\mathbb{Q}$} in: Yui, N., Lewis, J. D. (eds.): \emph{Calabi-Yau Varieties and Mirror %%@
Symmetry} (Toronto 2001), Fields Inst. Comm. 38, AMS (2003), 159-166.
\bibitem[HV]{HV} Hulek, K., Verrill, H.: \emph{On modularity of rigid and nonrigid Calabi-Yau %%@
varieties associated to the root lattice $A_4$}, arXiv: {\tt math.AG/0304169}
\bibitem[J] {J} Jones, J.: \emph{Tables of number fields with prescribed ramification}, {\tt %%@
http://math.la.asu.edu/$\sim$jj/numberfields}.
\bibitem[L] {L} Livn\'e, R.: \emph{Cubic exponential sums and Galois representations}, in: K. %%@
Ribet (Hrsg.), \emph{Current trends in arithmetical algebraic geometry} (Arcata 1985), Contemp. %%@
Math. 67, Amer. Math. Soc. (1987), 247-261.
%\bibitem[MP]{MP} Miranda, R., Persson, U.: \emph{On Extremal Rational Elliptic Surfaces}, Math. %%@
%Z. 193 (1986), 537-558.
\bibitem[S]{S} Sch\"utt, Matthias: \emph{New examples of modular rigid Calabi-Yau threefolds}, %%@
Collect. Math. 55,2 (2004), 219-228.
\bibitem[Sc1] {Sc1} Schoen, C.: \emph{On the geometry of a special determinantal hypersurface %%@
associated to the
 Mumford-Horrocks vector bundle}, J. reine angew. Math. 364 (1986), 85-111.
\bibitem[Sc2] {Sc2} Schoen, C.: \emph{On Fiber Products of Rational Elliptic Surfaces with %%@
Section}, Math. Z. 197 (1988), 177-199.
\bibitem[Se1] {Se1} Serre, J.-P.: \emph{R\'esum\'e des cours de 1984-1985}, Annuaire du Collège de %%@
France (1985), 85-90.
\bibitem[Se2] {Se2} Serre, J.-P.: \emph{Abelian $\ell$-Adic Representations and Elliptic Curves}, Research Notes in Mathematics, Vol. 7, A K Peters, Wellesley (1998).
\bibitem[St] {St} Stein, W.: \emph{Modular forms database}, {\tt http://modular.fas.harvard.edu}.
\bibitem[SY] {SY} Saito, M.-H., Yui, N.: \emph{The modularity conjecture for rigid Calabi-Yau %%@
threefolds over $\mathbb{Q}$}, Kyoto J. Math. 41, no. 2 (2001), 403-419.
\bibitem[W] {W} Wiles, A.: \emph{Modular elliptic curves and Fermat´s Last Theorem}, Ann. Math. %%@
(2) 141, No. 3 (1995), 443-551.
\bibitem[Y]{Y} Yui, Noriko: \emph{Update on the Modularity of Calabi-Yau Varieties}, in: Yui, N., %%@
Lewis, J. D. (eds.): \emph{Calabi-Yau Varieties and Mirror Symmetry} (Toronto 2001), Fields Inst. %%@
Comm. 38, AMS (2003), 307-362.

\end{thebibliography}}
\end{document}